\documentclass[Svgc]{Svgc}
\usepackage{amsmath,graphicx}
\usepackage[sort&compress,numbers]{natbib}
\journalname{Graphs and Combinatorics}

\newcommand{\thmlabel}[1]{\label{thm:#1}}
\newcommand{\thmref}[1]{Theorem~\ref{thm:#1}}
\newcommand{\lemlabel}[1]{\label{lem:#1}}
\newcommand{\lemref}[1]{Lemma~\ref{lem:#1}}

\newcommand{\eqnlabel}[1]{\label{eqn:#1}}
\newcommand{\eqnref}[1]{\eqref{eqn:#1}}
\newcommand{\Eqnref}[1]{Equation~\eqref{eqn:#1}}
\newcommand{\figlabel}[1]{\label{fig:#1}}
\newcommand{\figref}[1]{Figure~\ref{fig:#1}}
\newcommand{\seclabel}[1]{\label{sec:#1}}
\newcommand{\applabel}[1]{\label{app:#1}}
\newcommand{\secref}[1]{Section~\ref{sec:#1}}
\newcommand{\appref}[1]{Appendix~\ref{app:#1}}
\newcommand{\corlabel}[1]{\label{cor:#1}}
\newcommand{\corref}[1]{Corollary~\ref{cor:#1}}

\newcommand{\proplabel}[1]{\label{prop:#1}}
\newcommand{\propref}[1]{Proposition~\ref{prop:#1}}
\newcommand{\twopropref}[2]{Propositions~\ref{prop:#1} and ~\ref{prop:#2}}

\newcommand{\FLOOR}[1]{\ensuremath{\protect\left\lfloor#1\right\rfloor}}
\newcommand{\floor}[1]{\ensuremath{\protect\lfloor#1\rfloor}}
\newcommand{\Oh}[1]{\ensuremath{\protect\mathcal{O}(#1)}}
\newcommand{\bracket}[1]{\ensuremath{\protect\left(#1\right)}}
\newcommand{\SET}[1]{\ensuremath{\protect\left\{#1\right\}}}

\newcommand{\half}{\ensuremath{\protect\tfrac{1}{2}}}

\newcommand{\Figure}[4][htb]{
\begin{figure}[#1]
	\begin{center}#3\end{center}
	\caption{\figlabel{#2}#4}
\end{figure}}

\begin{document}
\title{On the Maximum Number of Cliques in a Graph}
\author{David R. Wood\thanks{Research supported by a Marie Curie Fellowship of the European Community under contract 023865, and by the projects MCYT-FEDER BFM2003-00368 and Gen.\ Cat 2001SGR00224.}}
\institute{Departament de Matem{\'a}tica Aplicada II, Universitat Polit{\`e}cnica de Catalunya, Barcelona, Spain \email{david.wood@upc.es}}

\maketitle

\begin{abstract} 
A \emph{clique} is a set of pairwise adjacent vertices in a graph. We determine the maximum number of cliques in a graph for the following graph classes:
(1) graphs with $n$ vertices and $m$ edges;
(2) graphs with $n$ vertices, $m$ edges, and maximum degree $\Delta$;
(3) $d$-degenerate graphs with $n$ vertices and $m$ edges;
(4) planar graphs with $n$ vertices and $m$ edges; and
(5) graphs with $n$ vertices and no $K_5$-minor or no $K_{3,3}$-minor.
For example, the maximum number of cliques in a planar graph with $n$ vertices is $8(n-2)$. 
\end{abstract}


\receive{June 23, 2006}

\section{Introduction}
\seclabel{Intro}

The typical question of extremal graph theory asks for the maximum number of edges in a graph in a certain family; see the surveys \citep{SimonSos-DM01, Simonovits97, Bollobas95, Simonovits83}. For example, a celebrated theorem of \citet{Turan41} states that the maximum number of edges in a graph with $n$ vertices and no $(k+1)$-clique is $\half(1-\frac{1}{k})n^2$. Here a \emph{clique} is a (possibly empty) set of pairwise adjacent vertices in a graph. For $k\geq0$, a \emph{$k$-clique} is a clique of cardinality $k$. Since an edge is nothing but a $2$-clique, it is natural to consider the maximum number of $\ell$-cliques in a graph. The following generalisation of Tur\'{a}n's Theorem, first proved by \citet{Zykov49}, has been rediscovered and itself generalised by several authors \citep{FisherRyan-DM92, HadNen77, PR00, NenKhad75, KN75, Had77, Roman-DM76, Erdos-CPM69, Erdos62, Sauer-JCTB71}.

\begin{theorem}[\citep{Zykov49}]
\thmlabel{Zykov}
For all integers $k\geq\ell\geq 0$, the maximum number of $\ell$-cliques in a graph with $n$ vertices and no $(k+1)$-clique is $\binom{k}{\ell}\bracket{\frac{n}{k}}^{\ell}$.
\end{theorem}

A simple inductive proof of \thmref{Zykov} is included in \appref{BoundedClique}. In this paper we determine the maximum number of cliques in a graph in each of the following classes:

\begin{itemize}

\item graphs with $n$ vertices and $m$ edges (\secref{General}),

\item graphs with $n$ vertices, $m$ edges, and maximum degree $\Delta$ (\secref{Degree}),

\item $d$-degenerate graphs with $n$ vertices and $m$ edges (\secref{Degenerate}),

\item planar graphs with $n$ vertices and $m$ edges (\secref{Planar}), and

\item graphs with $n$ vertices and no $K_5$-minor or no $K_{3,3}$-minor (\secref{K5Minor}).

\end{itemize}

We now review some related work from the literature. \citet{Eckhoff-DM04, Eckhoff-DM99} determined the maximum number of cliques in a graph with $m$ edges and no $(k+1)$-clique. Lower bounds on the number of cliques in a graph have also been obtained \citep{FisherNonis90, Moon-CMB65, Larman69, LovSim76, LovSim83, Fisher-JGT89, BolErdSze-DM75}. The number of cliques in a random graph has been studied \citep{OT-MS97, Schurger-PMH79, BolErd76}. Bounds on the number of cliques in a graph have recently been applied in the analysis of an algorithm for finding small separators \citep{ReedWood-EuroComb05} and in the enumeration of minor-closed families \citep{NSTW-JCTB06}.

\section{Preliminaries}
\seclabel{Prelim}

Every graph $G$ that we consider is undirected, finite, and simple. Let $V(G)$ and $E(G)$ be the vertex and edge sets of $G$. Let $\Delta(G)$ be the maximum degree of $G$. We say $G$ is a $(|V(G)|,|E(G)|)$-graph or a
$(|V(G)|,|E(G)|,\Delta(G))$-graph. 

Let $C(G)$ be the set of cliques in $G$. Let $c(G):=|C(G)|$. Let $C_k(G)$ be the set of $k$-cliques in $G$. Let $c_k(G):=|C_k(G)|$. Our aim is to prove bounds on $c(G)$ and $c_k(G)$. 

A clique is not necessarily maximal\footnote{\citet{MoonMoser} proved that the maximum number of \emph{maximal} cliques in a graph with $n$ vertices is approximately $3^{n/3}$; see \citep{GKSV05, SV05, Wilf-SJDM86, HL74, Spencer71, Erdos-IJM66, JC00, FHT-Networks93, FL91, Sato-DM86} for related results.}. In particular, $\emptyset$ is a clique of every graph, $\{v\}$ is a clique for each vertex $v$, and each edge is a clique. Thus every graph $G$ satisfies 
\begin{equation}
\eqnlabel{LowerBound}
c(G)\geq c_0(G)+c_1(G)+c_2(G)=1+|V(G)|+|E(G)|\enspace.
\end{equation}
A \emph{triangle} is a $3$-clique. \Eqnref{LowerBound} implies  that
\begin{equation}
\eqnlabel{TriangleFree}
c(G)=1+|V(G)|+|E(G)|\text{ if and only if $G$ is triangle-free.}
\end{equation}

Triangle-free graphs have the fewest cliques. Obviously the complete graph $K_n$ has the most cliques for a graph on $n$ vertices. In particular, $c(K_n)=2^n$ since every set of vertices in $K_n$ is a clique. 

Say $v$ is a vertex of a graph $G$. Let $G_v$ be the subgraph of $G$ induced by the neighbours of $v$. Observe that $X$ is a clique of $G$ containing $v$ if and only if $X=Y\cup\{v\}$ for some clique $Y$ of $G_v$. Thus the number of cliques of $G$ that contain $v$ is exactly $c(G_v)$. Every clique of $G$ either contains $v$ or is a clique of $G\setminus v$. Thus $C(G)=C(G\setminus v)\cup\{Y\cup\{v\}:Y\in C(G_v)\}$ and 
\begin{equation}
\eqnlabel{DeleteVertex}
c(G)=c(G\setminus v)+c(G_v)\leq c(G\setminus v)+2^{\deg(v)}\enspace.
\end{equation}

Let $G$ be a graph with induced subgraphs $G_1$, $G_2$ and $S$ such that $G=G_1\cup G_2$ and $G_1\cap G_2=S$. Then $G$ is obtained by \emph{pasting} $G_1$ and $G_2$ on $S$. Observe that $C(G)=C(G_1)\cup C(G_2)$ and $C(G_1)\cap C(G_2)=C(S)$. Thus
\begin{equation}
\eqnlabel{Pasting}
c(G)=c(G_1)+c(G_2)-c(S)\enspace.
\end{equation}

\begin{lemma}
\lemlabel{PastingPasting}
Let $G$ be an $(n,m)$-graph that is obtained by pasting $G_1$ and $G_2$ on $S$. Say $G_i$ has $n_i$ vertices and $m_i$ edges. Say $S$ has $n_S$ vertices and $m_S$ edges. If $c(G_i)\leq xn_i+ym_i+z$ and $c(S)\geq xn_s+ym_S+z$, then $c(G)\leq xn+ym+z$.
\end{lemma}

\begin{proof}
By \Eqnref{Pasting}, 
\begin{align*}
c(G)
& = c(G_1)+c(G_2)-c(S)\\
& \leq (xn_1+ym_1+z)+(xn_2+ym_2+z)-(xn_s+ym_s+z)\\
& = x(n_1+n_2-n_S)+y(m_1+m_2-m_S)+z\\
& = xn+ym+z\enspace.
\end{align*}
\qed\end{proof}

The following special case of \lemref{PastingPasting} will be useful.

\begin{corollary}
\corlabel{PastingClique}
Let $G$ be an $(n,m)$-graph that is obtained by pasting $G_1$ and $G_2$ on a $k$-clique. Say $G_i$ has $n_i$ vertices and $m_i$ edges. Assume that $c(G_i)\leq xn_i+ym_i+z$ and that $xk+y\binom{k}{2}+z\leq 2^k$. Then $c(G)\leq xn+ym+z$.\qed
\end{corollary}

\section{General Graphs}
\seclabel{General}

We now determine the maximum number of cliques in an $(n,m)$-graph.

\begin{theorem}
\thmlabel{GeneralGraphs}
Let $n$ and $m$ be non-negative integers such that $m\leq\binom{n}{2}$. Let $d$ and $\ell$ be the unique integers such that $m=\binom{d}{2}+\ell$ where $d\geq1$ and $0\leq\ell\leq d-1$. Then the maximum number of cliques in an $(n,m)$-graph equals $2^d+2^\ell+n-d-1$.
\end{theorem}

\begin{proof}
First we prove the lower bound. Let $V(G):=\{v_1,v_2,\dots,v_n\}$ and $E(G):=\{v_iv_j:1\leq i<j\leq d\}\cup\{v_iv_{d+1}:1\leq i\leq\ell\}$, as illustrated in \figref{Example}. Then $G$ has $\binom{d}{2}+\ell$ edges. Now $\{v_1,v_2,\dots,v_d\}$ is a clique, which contains $2^d$ cliques (including $\emptyset$). The neighbourhood of $v_{d+1}$ is an $\ell$-clique with $2^{\ell}$ cliques. Thus there are $2^\ell$ cliques that contain $v_{d+1}$. Finally $v_{d+2},v_{d+3},\dots,v_n$ are isolated vertices, which contribute $n-d-1$ cliques to $G$. In total, $G$ has $2^d+2^\ell+n-d-1$ cliques.

\Figure{Example}{\includegraphics{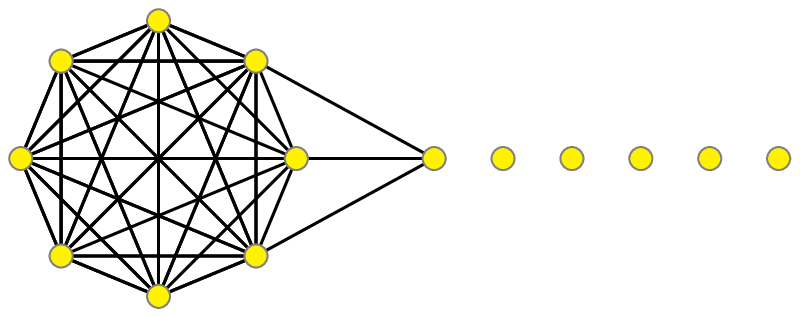}}{A $(14,31)$-graph with $269$ cliques ($d=8$ and $\ell=3$).}

Now we prove the upper bound. That is, every $(n,m)$-graph $G$ has at most $2^d+2^\ell+n-d-1$ cliques. We proceed by induction on $n+m$. For the base case, suppose that $m=0$. Then $d=1$, $\ell=0$, and $c(G)=n+1=2^d+2^\ell+n-d-1$. Now assume that $m\geq1$. Let $v$ be a vertex of minimum degree in $G$. Then $\deg(v)\leq d-1$, as otherwise every vertex has degree at least $d$, implying $m\geq\frac{dn}{2}\geq \frac{d(d+1)}{2}=\binom{d+1}{2}$, which contradicts the definition of $d$. By \Eqnref{DeleteVertex}, $c(G)\leq c(G\setminus v)+2^{\deg(v)}$. To apply induction to $G\setminus v$ (which has $n-1$ vertices and $m-\deg(v)$ edges) we distinguish two cases. 

First suppose that $\deg(v)\leq\ell$. Thus $m-\deg(v)=\binom{d}{2}+\ell-\deg(v)$. 
By induction,
$c(G)\leq 2^d+2^{\ell-\deg(v)}+n-1-d-1+2^{\deg(v)}$. 
Hence the result follows if
$2^d+2^{\ell-\deg(v)}+n-1-d-1+2^{\deg(v)}\leq2^d+2^\ell+n-d-1$. 
That is, 
$2^{\ell-\deg(v)}-1\leq(2^{\ell-\deg(v)}-1)2^{\deg(v)}$, which is true since $0\leq\deg(v)\leq\ell$.

Otherwise $\ell+1\leq\deg(v)\leq d-1$.  Thus $m-\deg(v)=\binom{d-1}{2}+d-1+\ell-\deg(v)$. 
By induction,
$c(G)\leq2^{d-1}+2^{d-1+\ell-\deg(v)}+n-1-d+2^{\deg(v)}$.
Hence the result follows if
$2^{d-1}+2^{d-1+\ell-\deg(v)}+n-1-d+2^{\deg(v)}\leq2^d+2^\ell+n-d-1$. 
That is, 
$2^\ell(2^{\deg(v)-\ell}-1) \leq 2^{d-1-\deg(v)+\ell}(2^{\deg(v)-\ell}-1)$. 
Since $\deg(v)\geq\ell+1$, we need
$2^\ell \leq 2^{d-1-\deg(v)+\ell}$, which is true since $\deg(v)\leq d-1$.
\qed\end{proof}

\section{Bounded Degree Graphs}
\seclabel{Degree}

We now determine the maximum number of cliques in an $(n,m,\Delta)$-graph.
\citet{West84} proved a related result.

\begin{theorem}
\thmlabel{Degree}
The number of cliques in an $(n,m,\Delta)$-graph $G$ is at most
\begin{equation*}
1+n+\bracket{\frac{2^{\Delta+1}-\Delta-2}{\binom{\Delta+1}{2}}}m\leq
1+\bracket{\frac{2^{\Delta+1}-1}{\Delta+1}}n
\enspace.
\end{equation*} 
\end{theorem}

\begin{proof}
$G$ has one $0$-clique and $n$ $1$-cliques. For $k\geq2$, each edge is in at most $\binom{\Delta-1}{k-2}$ $k$-cliques, and each $k$-clique contains $\binom{k}{2}$ edges. Thus $G$ has at most $m\binom{\Delta-1}{k-2}/\binom{k}{2}$ $k$-cliques. Thus the number of cliques (not counting $0$- and $1$-cliques) is at most
\begin{align*}
\sum_{k=2}^{\Delta+1}\frac{m\binom{\Delta-1}{k-2}}{\binom{k}{2}}
&=
m\sum_{k=2}^{\Delta+1}
\frac{2}{k(k-1)}\cdot\frac{(\Delta-1)!}{(k-2)!(\Delta-1-k+2)!}\\
&=
\frac{m}{\binom{\Delta+1}{2}}
\sum_{k\geq2}^{\Delta+1}\frac{2(\Delta-1)!\binom{\Delta+1}{2}}{k!(\Delta+1-k)!}\\
&=
\frac{m}{\binom{\Delta+1}{2}}
\sum_{k=2}^{\Delta+1}\frac{(\Delta+1)!}{k!(\Delta+1-k)!}\\
&=
\frac{m}{\binom{\Delta+1}{2}}
\bracket{
\bracket{\sum_{k=0}^{\Delta+1}\binom{\Delta+1}{k}}
-\frac{(\Delta+1)!}{1!(\Delta+1-1)!}
-\frac{(\Delta+1)!}{0!(\Delta+1-0)!}
}\\
&=
\frac{m}{\binom{\Delta+1}{2}}
\bracket{2^{\Delta+1}-\Delta-2}\enspace.
\end{align*}
The result follows since $m\leq\frac{\Delta n}{2}$. 
\qed\end{proof}

The bound in \thmref{Degree} is tight for many values of $m$.

\begin{proposition}
\proplabel{DegreeLowerBound}
For all $n$ and $m$ such that $m\leq\frac{\Delta n}{2}$ and $m\equiv0\pmod{\binom{\Delta+1}{2}}$, there is an $(n,m,\Delta)$-graph $G$ with
\begin{equation*}
c(G)=1+n+\bracket{\frac{2^{\Delta+1}-\Delta-2}{\binom{\Delta+1}{2}}}m\enspace.
\end{equation*}
\end{proposition}

\begin{proof}
Let $p:=m/\binom{\Delta+1}{2}$. Let $G$ consist of $p$ copies of $K_{\Delta+1}$, plus $n-p(\Delta+1)$ isolated vertices. Then $G$ is an $(n,m,\Delta)$-graph. Each copy of $K_{\Delta+1}$ contributes $2^{\Delta+1}-\Delta-2$ cliques with at least two vertices. Thus $G$ has $1+n+(2^{\Delta+1}-\Delta-2)p$ cliques.
\qed\end{proof}

\section{Degenerate Graphs}
\seclabel{Degenerate}

A graph $G$ is $d$-degenerate if every subgraph of $G$ has a vertex with degree at most $d$. The following simple result is well known; see \citep{ReedWood-EuroComb05, Eppstein-JGT93} for example. 

\begin{proposition}
\proplabel{Degen}
Every $d$-degenerate graph $G$ with $n\geq d$ vertices has at most $2^d(n-d+1)$ cliques.
\end{proposition}

\begin{proof}
We proceed by induction on $n$. If $n=d$ then $c(G)\leq2^d=2^d(n-d+1)$. Now assume that $n\geq d+1$. Let $v$ be a vertex of $G$ with $\deg(v)\leq d$. 
By \Eqnref{DeleteVertex}, $c(G)\leq c(G\setminus v)+2^{\deg(v)}$. Now $G\setminus v$ is $d$-degenerate since it is a subgraph of $G$. Moreover, $G\setminus v$ has at least $d$ vertices. By induction, $c(G\setminus v)\leq 2^d(n-1-d+1)$. Thus $c(G)\leq2^d(n-1-d+1)+2^d=2^d(n-d+1)$.
\qed\end{proof}

The bound in \propref{Degen} is tight.

\begin{proposition}
\proplabel{DegenLowerBound}
For all $n\geq d$, there is a $d$-degenerate graph $G_n$ with $n$ vertices and exactly $2^d(n-d+1)$ cliques \textup{(}and with a $d$-clique\textup{)}.
\end{proposition}

\begin{proof}
Let $G_d$ be the complete graph $K_d$. Then $G_d$ has the desired properties. For $n\geq d+1$, let $G_n$ be the graph obtained by adding one new vertex $v$ adjacent to every vertex in some $d$-clique in $G_{n-1}$. Then $G_n$ is $d$-degenerate and contains a $d$-clique. ($G_n$ is a chordal graph called a \emph{$d$-tree}; see \citep{Bodlaender-TCS98}.)\ By \Eqnref{DeleteVertex}, $c(G_n)=c(G_{n-1})+2^{\deg(v)}=2^d(n-1-d+1)+2^d=2^d(n-d+1)$.
\qed\end{proof}


\propref{Degen} can be made sensitive to the number of edges as follows. 

\begin{theorem}
\thmlabel{SensitiveDegen}
For all $d\geq1$, every $d$-degenerate graph $G$ with $n$ vertices and $m\geq\binom{d}{2}$ edges has at most 
\begin{equation*}
n+\frac{(2^d-1)m}{d}-\frac{(d-3)2^d+d+1}{2}
\end{equation*} 
cliques.
\end{theorem}

\begin{proof}
We proceed by induction on $n+m$. For the base case, suppose that $m=\binom{d}{2}+\ell$ where $d\geq1$ and $0\leq\ell\leq d-1$. Thus $c(G)\leq 2^d+2^\ell+n-d-1$ by \thmref{GeneralGraphs}, and the result follows if 
\begin{equation*}
2^d+2^\ell+n-d-1\leq n+\frac{(2^d-1)m}{d}-\frac{(d-3)2^d+d+1}{2}\enspace.
\end{equation*} 
That is, 
$d(2^\ell-1)\leq\ell(2^d-1)$, which we prove in \lemref{Ineq} below.

Now assume that $m\geq\binom{d+1}{2}$. Now $G$ has a vertex $v$ with $\deg(v)\leq d$. By \Eqnref{DeleteVertex}, $c(G)\leq c(G\setminus v)+2^{\deg(v)}$. The graph $G\setminus v$ has $m-\deg(v)\geq\binom{d}{2}$ edges, and is $d$-degenerate since it is a subgraph of $G$. By induction,  
\begin{equation*}c(G\setminus v)\leq n-1+\frac{(2^d-1)(m-\deg(v))}{d}-\frac{(d-3)2^d+d+1}{2}\enspace.
\end{equation*}
Thus the result follows if
\begin{equation*}-1+\frac{(2^d-1)(m-\deg(v))}{d}+2^{\deg(v)}\leq \frac{(2^d-1)m}{d}\enspace.\end{equation*}
That is, $d(2^{\deg(v)}-1)\leq (2^d-1)\deg(v)$, which holds by \lemref{Ineq} below.
\qed\end{proof}

\begin{lemma}
\lemlabel{Ineq}
$d(2^\ell-1)\leq \ell(2^d-1)$ for all integers $d\geq\ell\geq0$.
\end{lemma}

\begin{proof}
The case $\ell=0$ is trivial. Now assume that $\ell\geq1$.
We proceed by induction on $d$. The base case $d=\ell$ is trivial. 
Assume that $d\geq\ell+1\geq2$ and by induction,
\begin{equation}
\eqnlabel{AAA}
(d-1)(2^\ell-1)\leq\ell(2^{d-1}-1).
\end{equation} 
Since $d\geq2$, 
\begin{equation}
\eqnlabel{BBB}
\frac{d}{d-1}\leq 2<2+\frac{1}{2^{d-1}-1}=\frac{2^d-1}{2^{d-1}-1}.
\end{equation}
Equations \eqnref{AAA} and \eqnref{BBB} imply that
\begin{equation*}
(d-1)(2^\ell-1)\cdot\frac{d}{d-1}<\ell(2^{d-1}-1)\cdot\frac{2^d-1}{2^{d-1}-1}.
\end{equation*}
That is, $d(2^\ell-1)<\ell(2^d-1)$, as desired.
\qed\end{proof}


Note that a $d$-degenerate $n$-vertex graph has at most $dn-\binom{d+1}{2}$ edges, and \thmref{SensitiveDegen} with $m=dn-\binom{d+1}{2}$ is equivalent to \propref{Degen}.


The bound in \thmref{SensitiveDegen} is tight for many values of $m$.

\begin{proposition}
\proplabel{SensitiveDegenLowerBound}
Let $d\geq1$. For all $n$ and $m$ such that
$\binom{d}{2}\leq m\leq dn-\binom{d+1}{2}$ and
\begin{equation*}
m\bmod{d}=
\begin{cases}
0			& \text{if }d\text{ is odd} \\
\frac{d}{2}	& \text{if }d\text{ is even ,}
\end{cases}
\end{equation*}
there is a $d$-degenerate $(n,m)$-graph $G$ with 
\begin{equation*}
c(G)=n+\frac{(2^d-1)m}{d}-\frac{(d-3)2^d+d+1}{2}\enspace.
\end{equation*}
\end{proposition}

\begin{proof}
Let $n':=\frac{m}{d}+\half(d+1)$. Then $n'$ is an integer and $d\leq n'\leq n$. Let $G$ consist of a $d$-degenerate $n'$-vertex graph with $2^d(n'-d+1)$ cliques (from \propref{DegenLowerBound}), plus $n-n'$ isolated vertices. Then $G$ has $m$ edges and $c(G)=2^d(n'-d+1)+n-n'=n+(2^d-1)\frac{m}{d}-\half((d-3)2^d+d+1)$. 
\qed\end{proof}


A graph is $1$-degenerate if and only if it is a forest. Thus 
\thmref{SensitiveDegen} with $d=1$ implies that every forest has at most $n+m-1$ cliques, which also follows from \Eqnref{TriangleFree}. In particular, $c(T)=2n$ for every $n$-vertex tree $T$. 

\thmref{SensitiveDegen} with $d=2$ implies that every $2$-degenerate graph has at most $n+\half(3m+1)$ cliques. Outerplanar graphs are $2$-degenerate. The construction in \twopropref{DegenLowerBound}{SensitiveDegenLowerBound} can produce outerplanar graphs. (Add each new vertex adjacent to two consecutive vertices on the outerface.)\ Thus this bound is tight for outerplanar graphs. 

\section{Planar Graphs}
\seclabel{Planar}

\citet{PY-IPL81} and \citet{Storch-GECCO06} proved that every $n$-vertex planar graph has \Oh{n} cliques; see \citep{Eppstein-JGT93} for a more general result. The proof is based on the corollary of Euler's Formula that planar graphs are $5$-degenerate. By \thmref{SensitiveDegen}, if $G$ is a planar $(n,m)$-graph with $m\geq 10$, then $c(G)<n+\frac{31}{5}m<\frac{98}{5}n$. We now prove that the bound for $3$-degenerate graphs in \thmref{SensitiveDegen} also holds for planar graphs.

\begin{theorem}
\thmlabel{PlanarGraphs}
Every planar $(n,m)$-graph $G$ with $m\geq3$ has at most $n+\frac{7}{3}m-2$ cliques.
\end{theorem}

\begin{proof}
We proceed by induction on $n+m$. The result is easily verified if $m=3$. 

Suppose that $G$ has a separating triangle $T$. Thus $G$ is obtained by pasting two induced subgraphs $G_1$ and $G_2$ on $T$. Say $G_i$ has $n_i$ vertices and $m_i$ edges. Then $m_i\geq3$ since $T\subset G_i$. By induction, $c(G_i)\leq n_i+\frac{7}{3}m_i-2$. By \corref{PastingClique} with $k=3$, $x=1$, $y=\frac{7}{3}$ and $z=-2$, we have $c(G)\leq n+\frac{7}{3}m-2$ (since $1\cdot3+\frac{7}{3}\binom{3}{2}-2=2^3$). Now assume that $G$ has no separating triangle. 

Let $v$ be a vertex of $G$. We have $c(G)=c(G\setminus v)+c(G_v)$ by \Eqnref{DeleteVertex}. The graph $G\setminus v$ has $m-\deg(v)$ edges. Suppose that $m-\deg(v)\leq 2$. (Then we cannot apply induction to $G\setminus v$.)\ Then $G$ has no $4$-clique and at most two triangles. If $G$ has at most one triangle, then $c(G)\leq 1+n+m+1\leq n+\frac{7}{3}m-2$ since $m\geq3$. Otherwise $G$ has two triangles, and $c(G)\leq 1+n+m+2<n+\frac{7}{3}m-2$ since $m\geq5$.

Now assume that $m-\deg(v)\geq3$. By \eqnref{DeleteVertex}, applying induction to $G\setminus v$, 
\begin{equation*}
c(G)=c(G\setminus v)+c(G_v)
\leq(n-1)+\tfrac{7}{3}(m-\deg(v))-2+c(G_v)\enspace.
\end{equation*}
Fix a plane embedding of $G$. If $uw$ is an edge of $G_v$, then the edges $vu$ and $vw$ are consecutive in the circular ordering of edges incident to $v$ defined by the embedding (as otherwise $G$ would contain a separating triangle). Thus $\Delta(G_v)\leq 2$ and $c(G_v)\leq1+\frac{7}{3}\deg(v)$ by \thmref{Degree}. Hence
\begin{equation*}
c(G)\leq(n-1)+(\tfrac{7}{3}(m-\deg(v))-2)+(1+\tfrac{7}{3}\deg(v))
=n+\tfrac{7}{3}m-2\enspace.
\end{equation*}
\qed\end{proof}

If $n\geq3$ in \thmref{PlanarGraphs} then $m\leq3(n-2)$ by Euler's Formula. Thus we have the following corollary.

\begin{corollary}
\corlabel{MaximalPlanarGraphs}
Every planar graph with $n\geq3$ vertices has at most $8(n-2)$ cliques.\qed
\end{corollary}

We now prove bounds on the number of $3$- and $4$-cliques in a planar graph.

\begin{proposition}
\proplabel{PlanarTrianglesSquares}
For every planar graph $G$ with $n\geq3$ vertices, $c_3(G)\leq 3n-8$ and  $c_4(G)\leq n-3$. 
\end{proposition}

\begin{proof}
We proceed by induction on $n$. The result is trivial if $n\leq 4$. Now assume that $n\geq5$. First suppose that $G$ has no separating triangle. Then $c_4(G)=0$, and every triangle of $G$ is a face. By Euler's Formula, $c_3(G)\leq 2n-4<3n-8$ faces. Now suppose that $G$ has a separating triangle $T$. Thus $G$ is obtained by pasting two induced subgraphs $G_1$ and $G_2$ on $T$. Say $G_i$ has $n_i$ vertices. Then $n_i\geq3$ since $T\subset G_i$. By induction, $c_3(G_i)\leq 3n_i-8$ and $c_4(G_i)\leq n_i-3$. Every clique of $G$ is a clique of $G_1$ or $G_2$. Thus $c_4(G)=c_4(G_1)+c_4(G_2)\leq n_1-3+n_2-3=n-3$. Moreover, $T$ is a triangle in both $G_1$ and $G_2$. Thus 
$c_3(G)\leq (3n_1-8)+(3n_2-8)-1=3(n_1+n_2)-17=3(n+3)-17=3n-8$.
\qed\end{proof}

Note that \propref{PlanarTrianglesSquares} and Euler's Formula (which implies $c_2(G)\leq 3n-6$) reprove \corref{MaximalPlanarGraphs}, since $1+n+3(n-2)+(3n-8)+(n-3)=8(n-2)$. 

We now show that all our bounds for planar graphs are tight.

\begin{proposition}
\proplabel{MaximalPlanarLowerBound}
For all $n\geq3$ there is a maximal planar $n$-vertex graph $G_n$ with
$c_2(G_n)=3(n-2)$, $c_3(G_n)=3n-8$, $c_4(G_n)=n-3$, and $c(G_n)=8(n-2)$.
\end{proposition}

\begin{proof}
Let $G_3:=K_3$. Then $c_2(G_3)=3$, $c_3(G_3)=1$, $c_4(G_3)=0$, and $c(G_3)=8$. Say $G_{n-1}$ is a maximal planar $(n-1)$-vertex graph with 
$c_2(G_{n-1})=3(n-3)$, $c_3(G_{n-1})=3n-11$, $c_4(G_{n-1})=n-4$, and $c(G_n)=8(n-3)$. Let $G_n$ be the maximal planar $n$-vertex graph obtained by adding one new vertex $v$ adjacent to each vertex of some face of $G_{n-1}$, as illustrated in \figref{Stellated}. Then
$c_2(G_n)=c_2(G_{n-1})+3=3(n-2)$, 
$c_3(G_n)=c_3(G_{n-1})+3=3n-8$, 
$c_4(G_n)=c_4(G_{n-1})+1=n-3$, and $c(G_n)=c(G_{n-1})+c(G_n(v))=8(n-3)+8=8(n-2)$. (Note that $G_n$ is also an example of a $3$-degenerate graph with the maximum number of cliques; see \propref{DegenLowerBound}.)\ 
\qed\end{proof}

\Figure{Stellated}{\includegraphics[width=\textwidth]{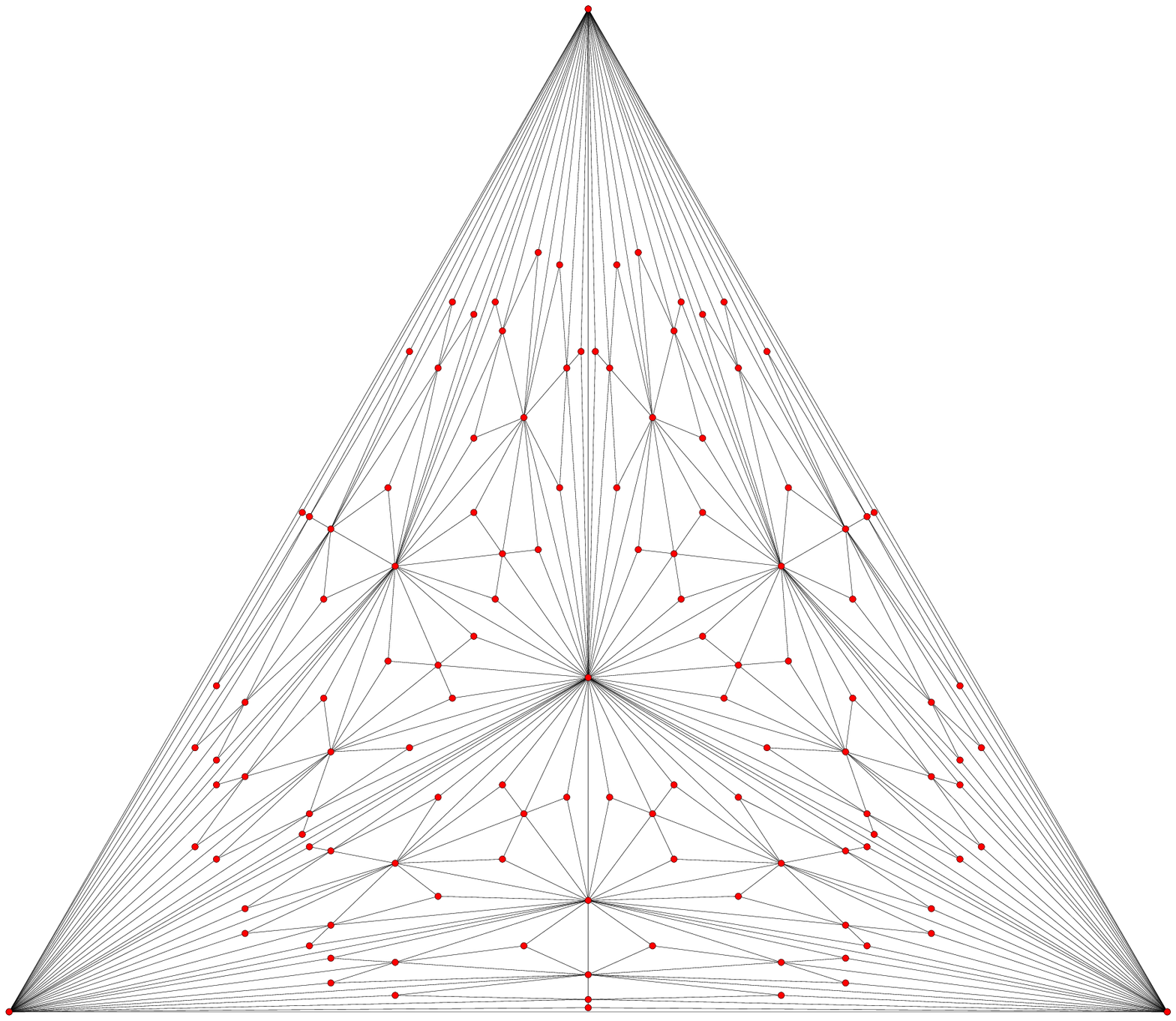}}{A planar graph with $124$ vertices, $366$ edges, $364$ triangles, $121$ $4$-cliques, and $976$ cliques. It is obtained by repeatedly adding one degree-$3$ vertex inside each internal face (starting from $K_3$).}

\begin{proposition}
\proplabel{PlanarLowerBound}
For all $n\geq3$ and $m\in\{3,6,\dots,3n-6\}$, there is a planar $(n,m)$-graph $G$ with $c(G)=n+\frac{7}{3}m-2$.
\end{proposition}

\begin{proof}
Let $n':=\frac{m}{3}+2$. Let $G$ consist of a maximal planar graph on $n'$ vertices with $8(n'-2)$ cliques (from \propref{MaximalPlanarLowerBound}), plus $n-n'$ isolated vertices. Then $G$ has $n$ vertices and $m$ edges, and $c(G)=8(n'-2)+n-n'=n+7n'-16=n+7(\frac{m}{3}+2)-16=n+\frac{7}{3}m-2$.
\qed\end{proof}



\section{Graphs with no $K_5$-Minor}
\seclabel{K5Minor}

A graph $H$ is a \emph{minor} of a graph $G$ if $H$ can be obtained from a subgraph of $G$ by contracting edges. The graphs with no $K_3$-minor are the forests, which have at most $2n$ cliques, and this bound is tight. The graphs with no $K_4$-minor (called \emph{series-parallel}) are $2$-degenerate, and thus have at most $4(n-1)$ cliques, and this bound is tight. The Kuratowski-Wagner Theorem characterises planar graphs as those with no $K_5$-minor and no $K_{3,3}$-minor. We now extend \corref{MaximalPlanarGraphs} for graphs with no $K_5$-minor (but possibly a $K_{3,3}$-minor).

\begin{theorem}
\thmlabel{NoK5Minor}
Every graph $G$ with $n\geq3$ vertices and no $K_5$-minor has at most $8(n-2)$ cliques.
\end{theorem}

\begin{proof}
Let $V_8$ be the graph obtained from the $8$-cycle by adding an edge between each pair of antipodal vertices; see \figref{W}. Let $G$ be a minimum counterexample to the theorem. We can assume that $G$ is edge-maximal with no $K_5$-minor. \citet{Wagner37} proved that (a) $G$ is a maximal planar graph, (b) $G=V_8$, or (c) $G$ is obtained by pasting two smaller graphs (that are thus not counterexamples), each with no $K_5$-minor, on an edge or a triangle $T$. In case (a) the result is \corref{MaximalPlanarGraphs}. In case (b), since $V_8$ is triangle-free, $c(V_8)=1+|V(V_8)|+|E(V_8)|=21<8(|V(V_8)|-2)$ by \Eqnref{TriangleFree}. In case (c), if $T$ is an edge, we have $c(G)\leq 8(n-2)$ by \corref{PastingClique} with $k=2$, $x=8$, $y=0$ and $z=-16$  (since $8\cdot2+0-16<2^2$). In case (c), if $T$ is a triangle, we have $c(G)\leq 8(n-2)$ by \corref{PastingClique} with $k=3$, $x=8$, $y=0$ and $z=-16$ (since $8\cdot3+0-16=2^3$). 
\qed\end{proof}

\Figure{W}{\includegraphics{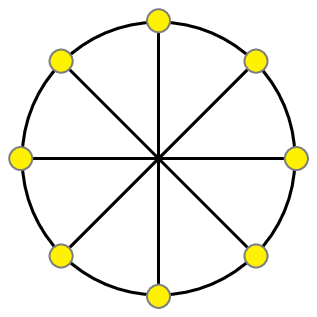}}{The graph $V_8$.}

A similar result is obtained for graphs with no $K_{3,3}$-minor.

\begin{theorem}
\thmlabel{NoK33Minor}
Every graph $G$ with $n\geq3$ vertices and no $K_{3,3}$-minor has at most $\frac{4}{3}(7n-11)$ cliques. Conversely, for all $n\equiv 2\pmod{3}$ with $n\geq5$ there is an $n$-vertex graph with no $K_{3,3}$-minor and $c(G)=\frac{4}{3}(7n-11)$. 
\end{theorem}

\begin{proof}
Let $G$ be a minimum counterexample. We can assume that $G$ is edge-maximal with no $K_{3,3}$-minor. \citet{Wagner37} proved that (a) $G$ is a maximal planar graph, (b) $G=K_5$, or (c) $G$ is obtained by pasting two smaller graphs (that are thus not counterexamples), each with no $K_{3,3}$-minor, on an edge. In case (a) the result follows from  \corref{MaximalPlanarGraphs} since $8n-16<\frac{4}{3}(7n-11)$. In case (b), $c(K_5)=32=\frac{4}{3}(7\cdot 5-11)$. In case (c), we have $c(G)\leq \frac{4}{3}(7n-11)$ by \corref{PastingClique} with $k=2$, $x=\frac{28}{3}$, $y=0$ and $z=-\frac{44}{3}$  (since $\frac{28}{3}\cdot2+0-\frac{44}{3}=2^2$). By the same analysis, the graph obtained from $K_5$ by repeatedly pasting copies of $K_5$ on an edge has no $K_{3,3}$-minor and $\frac{4}{3}(7n-11)$ cliques.
\qed\end{proof}


We finish with an open problem: What is the maximum number of cliques in an
$n$-vertex graph $G$ with no $K_t$-minor? \citet{Kostochka82} and \citet{Thomason84} independently proved that $G$ is \Oh{t\sqrt{\log t}}-degenerate\footnote{Moreover, this bound is best possible; \citet{Thomason01} even determined the asymptotic constant.}. Thus \propref{Degen} implies that $G$ has at most $2^{\Oh{t\sqrt{\log t}}}n$ cliques; similar bounds can be found in \citep{NSTW-JCTB06, ReedWood-EuroComb05}. It is unknown whether this bound can be improved to $c^tn$ for some constant $c$ (possibly for sufficiently large $n$). 

We have proved that $c(G)\leq 2^{t-2}(n-t+3)$ whenever $t\leq5$. Moreover, the graph $G$ in \propref{DegenLowerBound} (with $t=d+2$) has no $K_t$-minor and $c(G)=2^{t-2}(n-t+3)$. However, for large values of $t$ this upper bound does not hold for the complete $k$-partite graph $K_{2,2,\dots,2}$. By  \thmref{HadwigerCompletePartite} in \appref{Hadwiger}, the maximum order of a clique minor in $K_{2,2,\dots,2}$ is $\floor{\frac{3}{2}k}$. But by \propref{CompletePartite}, $c(K_{2,2,\dots,2})=3^k>2^{\floor{3k/2}-1}(2k-\floor{\frac32k}+2)$ for all $k\geq 42$.

\begin{acknowledgement}
Thanks to a referee for pointing out reference \citep{Stiebitz-DM92}.
\end{acknowledgement}


\def\soft#1{\leavevmode\setbox0=\hbox{h}\dimen7=\ht0\advance \dimen7
  by-1ex\relax\if t#1\relax\rlap{\raise.6\dimen7
  \hbox{\kern.3ex\char'47}}#1\relax\else\if T#1\relax
  \rlap{\raise.5\dimen7\hbox{\kern1.3ex\char'47}}#1\relax \else\if
  d#1\relax\rlap{\raise.5\dimen7\hbox{\kern.9ex \char'47}}#1\relax\else\if
  D#1\relax\rlap{\raise.5\dimen7 \hbox{\kern1.4ex\char'47}}#1\relax\else\if
  l#1\relax \rlap{\raise.5\dimen7\hbox{\kern.4ex\char'47}}#1\relax \else\if
  L#1\relax\rlap{\raise.5\dimen7\hbox{\kern.7ex
  \char'47}}#1\relax\else\message{accent \string\soft \space #1 not
  defined!}#1\relax\fi\fi\fi\fi\fi\fi} \def\cprime{$'$}
  \def\Dbar{\leavevmode\lower.6ex\hbox to 0pt{\hskip-.23ex \accent"16\hss}D}

\appendix
\section{Graphs with Bounded Cliques}
\applabel{BoundedClique}


In this appendix we give a simple inductive proof of \thmref{Zykov}. 

\begin{proposition}
\proplabel{BoundedCliqueClique}
For all integers $k\geq\ell\geq 0$, every graph $G$ with $n\geq\ell$ vertices and no $(k+1)$-clique has at most $\binom{k}{\ell}\bracket{\frac{n}{k}}^{\ell}$ $\ell$-cliques.
\end{proposition}

\begin{proof}
We proceed by induction on $n$. For the base case, suppose that $n\leq k$. Trivially $c_\ell(G)\leq\binom{n}{\ell}$, which is at most $\binom{k}{\ell}\bracket{\frac{n}{k}}^\ell$ by \lemref{Claim} below. Now assume that the result holds for graphs with less than $n$ vertices, and $n>k$. Let $G$ be a graph with $n$ vertices, no $(k+1)$-clique, and with $c_\ell(G)$ maximum. We can add edges to $G$ until it contains a $k$-clique $X$. Every $\ell$-clique of $G$ is the union of some $i$-clique of $G\setminus X$ and some $(\ell-i)$-clique of $G[X]$, for some $0\leq i\leq\ell$. Moreover, the vertices in each $i$-clique of $G\setminus X$ have at most $k-i$ common neighbours in $X$ (since $X$ is a clique and $G$ has no $(k+1)$-clique). Thus from each $i$-clique of $G\setminus X$, we obtain at most $\binom{k-i}{\ell-i}$ $\ell$-cliques of $G$. By induction, $c_i(G\setminus X)\leq\binom{k}{i}\bracket{\frac{n-k}{k}}^i$. Thus  
\begin{equation*}c_\ell(G)\leq
\sum_{i=0}^{\ell}\binom{k}{i}\bracket{\frac{n-k}{k}}^i\binom{k-i}{\ell-i}
=\binom{k}{\ell}\sum_{i=0}^{\ell}\binom{\ell}{i}\bracket{\frac{n}{k}-1}^i
=\binom{k}{\ell}\bracket{\frac{n}{k}}^{\ell}\enspace,
\end{equation*}
by the binomial theorem\footnote{Twice we use that $x^t=\sum_{j=0}^{t}\binom{t}{j}(x-1)^j$ for all real $x$.}.
\qed\end{proof}

\begin{lemma}
\lemlabel{Claim}
$\binom{n}{\ell}k^\ell\leq\binom{k}{\ell}n^\ell$ for all integers $k\geq n\geq\ell\geq0$.
\end{lemma}

\begin{proof}
We proceed by induction on $\ell$. The claim is trivial with $\ell=0$. Now assume that $\ell\geq1$. Thus $k-n\leq\ell(k-n)$, implying $kn+k-n\leq kn+\ell(k-n)$. That is, $k(n-\ell+1)\leq n(k-\ell+1)$. By induction,
\begin{equation*}
\binom{n}{\ell-1}k^{\ell-1}\cdot k(n-\ell+1)\leq
\binom{k}{\ell-1}n^{\ell-1}\cdot n(k-\ell+1)\enspace.
\end{equation*}
That is,
\begin{equation*}
\frac{n!\,k^\ell(n-\ell+1)}{(n-\ell+1)!\,(\ell-1)!}\leq
\frac{k!\,n^\ell(k-\ell+1)}{(k-\ell+1)!\,(\ell-1)!}\enspace.
\end{equation*}
Hence
\begin{equation*}
\frac{n!\,k^\ell}{(n-\ell)!\,\ell!}\leq
\frac{k!\,n^\ell}{(k-\ell)!\,\ell!}\enspace,
\end{equation*}
as desired.
\qed\end{proof}


\begin{proposition}
\proplabel{BoundedClique}
Every graph $G$ with $n$ vertices and no $(k+1)$-clique has at most $\bracket{\frac{n}{k}+1}^k$ cliques.
\end{proposition}

\begin{proof}
By \propref{BoundedCliqueClique} and the binomial theorem, 
\begin{equation*}
c(G)\leq 
\sum_{\ell=0}^k\binom{k}{\ell}\bracket{\frac{n}{k}}^{\ell}
=\bracket{\frac{n}{k}+1}^k\enspace.\end{equation*}
\qed\end{proof}

We now prove that \twopropref{BoundedCliqueClique}{BoundedClique} are tight.

\begin{proposition}
\proplabel{CompletePartite}
For every complete $k$-partite graph $G=K_{n_1,n_2,\dots,n_k}$, 
\begin{equation*}
c(G)\;=\;\prod_{i=1}^k(n_i+1)\enspace.
\end{equation*}
In particular, if every $n_i=\frac{n}{k}$ then 
$c(G)=(\frac{n}{k}+1)^k$ and $c_\ell(G)=\binom{k}{\ell}(\frac{n}{k})^\ell$ whenever $0\leq\ell\leq k$.
\end{proposition}

\begin{proof}
Every clique consists of at most one vertex from each of the $k$ colour classes. There are $n_i+1$ ways to choose at most one vertex from the $i$-th colour class. Thus $c(G)=\prod_i(n_i+1)$. (This result can also be proved using 
\Eqnref{DeleteVertex}.)\ Now assume that every $n_i=\frac{n}{k}$. Every $\ell$-clique consists of exactly one vertex from each of $\ell$ colour classes. There are $\binom{k}{\ell}$ ways to choose $\ell$ colour classes and $\frac{n}{k}$ ways to choose exactly one vertex from each colour class. Each combination gives a distinct $\ell$-clique. The result follows.
\qed\end{proof}

It is interesting to note that the extremal examples in \propref{DegreeLowerBound} for graphs of bounded degree (disjoint copies of cliques) are the complements of the extremal examples in \propref{CompletePartite} for graphs with bounded cliques (complete multipartite graphs).

\section{Clique Minors in a Complete Multipartite Graph}
\applabel{Hadwiger}

The \emph{Hadwiger number} of a graph $G$, denoted by $\eta(G)$, is the maximum order of a clique minor in $G$. \citet{Stiebitz-DM92} proved that $\eta(G)\leq\half(n+k)$ for every $n$-vertex graph $G$ with no $(k+1)$-clique. We now prove that this bound is tight for every complete $k$-partite graph if the largest colour class is not too large.

\begin{theorem}
\thmlabel{HadwigerCompletePartite}
Let $G$ be a complete $k$-partite graph on $n$ vertices with $n'$ vertices in the largest colour class. Then 
$\eta(G)=\min\SET{\half(n+k),n-n'+1}$.
\end{theorem}

The proof of \thmref{HadwigerCompletePartite} is based on the following lemma.

\begin{lemma}
\lemlabel{HadwigerCompletePartiteMatching}
Let $G$ be the complete $k$-partite graph $K_{n_1,n_2,\dots,n_k}$ with each $n_i\geq1$. Then $\eta(G)$ equals $k$ plus the size of the largest matching in $G':= K_{n_1-1,n_2-1,\dots,n_k-1}$.
\end{lemma}

\begin{proof}
Consider $G'$ to be a subgraph of $G$, so that $S:= V(G)\setminus V(G')$ is a $k$-clique of $G$. Let $M$ be a matching of $G'$. If $v$ is a vertex and $e$ is an edge of $G'$, then $v$ is adjacent to at least one endpoint of $e$. Thus every vertex in $S$ is adjacent to at least one endpoint of every edge in $M$, and for all edges $e$ and $f$ in $M$, at least one endpoint of $e$ is adjacent to at least one endpoint of $f$. Thus by contracting each edge of $M$ within $G$, we obtain a $K_{k+|M|}$-minor in $G$. 

Now suppose that $K_t$ is a minor of $G$ with $t$ maximum. Then $G$ has disjoint vertex sets $X_1,X_2,\dots,X_t$, such that each $X_i$ induces a connected subgraph of $G$, and for all $i\ne j$, some vertex in $X_i$ is adjacent to some vertex in $X_j$.

Suppose that some $X_i$ contains two vertices $v$ and $w$ in the same colour class of $G$. Since $v$ and $w$ have the same neighbourhood, we can delete $w$ from $X_i$ and still have a $K_t$-minor. Now assume that the vertices in each set $X_i$ are from distinct colour classes. 

Suppose that some $X_i$ contains at least three vertices $u,v,w$. Since the neighbourhood of $u$ is contained in the union of the neighbourhoods of $v$ and $w$, we can delete $u$ from $X_i$ and still have a $K_t$-minor. Now assume that each set $X_i$ has cardinality $1$ or $2$. 

Suppose that for some colour class $\ell$, no set $X_i$ contains a vertex coloured $\ell$. Then $X_1,\dots,X_t$ along with a set consisting of one vertex coloured $\ell$ forms a $K_{t+1}$-minor, which is a contradiction. Now assume that for every colour class $\ell$, there is some set $X_i$ that contains a vertex coloured $\ell$. 

Suppose that for some colour class $\ell$, every set $X_i$ that contains some vertex coloured $\ell$ has cardinality $2$. Let $X_i=\{v,w\}$ be such a set, where $v$ is coloured $\ell$. Thus $v$ is adjacent to some vertex in every set $X_j$. Thus we can delete $w$ from $X_i$ and still have a $K_t$-minor. Now assume that for each colour class $\ell$, some set $X_i$ consists of one vertex coloured $\ell$. No two singleton sets $X_i$ and $X_j$ contain vertices of the same colour. Thus there are $k$ singleton sets $X_i$, one for each colour class. The remaining sets $X_i$ thus form a matching in $G'$. 
\qed\end{proof}

\noindent\emph{Proof of \thmref{HadwigerCompletePartite}.} \citet{Sitton-EJUM} proved that the size of the largest matching in a complete multipartite graph on $n$ vertices with $n'$ vertices in the largest colour class is $\min\SET{\FLOOR{\frac{n}{2}},n-n'}$. Applying this result to the graph $G'$ in \lemref{HadwigerCompletePartiteMatching},
\begin{equation*}\eta(G)
\;=\;k+\min\SET{\half(n-k),(n-k)-(n'-1)}
\;=\;\min\SET{\half(n+k),n-n'+1}.
\end{equation*}
\qed

\end{document}